\documentclass[11pt]{article}
\usepackage{amsmath,amssymb,latexsym, color}
\usepackage{theorem}
\newtheorem{theorem}{Theorem}
\newtheorem{lemma}{Lemma}

\newtheorem{corollary}{Corollary}
\newtheorem{proposition}{Proposition}
{\theorembodyfont{\rmfamily} }

\title{A new approach to the asymptotics  for Sobolev  orthogonal polynomials}
\author{M. Alfaro$^{a,}$\thanks{Partially supported by MICINN of Spain under Grant
MTM2009-12740-C03-03, FEDER funds (EU), and the DGA project E-64
(Spain)}\,\,, J. J. Moreno--Balc\'{a}zar$^{b,}$\thanks{Partially
supported by MICINN of Spain under Grant MTM2008--06689--C02--01
and Junta de Andaluc\'{\i}a (FQM229 and excellence project
P06-FQM-1735).}\,\,, A. Pe\~{n}a$^{a,*}$, M. L. Rezola$^{a,*}$
 \\ \\ $^{a}$ {\small Departamento de Matem\'{a}ticas and IUMA.}\\ {Universidad
de Zaragoza (Spain).} \\ $^{b}$ {\small Departamento de
Estad\'{\i}stica y Matem\'{a}tica Aplicada and ICI.}\\{Universidad de Almer\'{\i}a (Spain).}}

\date{}

\begin{document}

\maketitle

\begin{abstract}
In this paper we deal with polynomials orthogonal with respect to
an inner product involving derivatives, that is, a Sobolev inner
product. Indeed, we consider Sobolev type polynomials which are
orthogonal with respect to $$(f,g)=\int fg d\mu +\sum_{i=0}^r M_i
f^{(i)}(0) g^{(i)}(0), \quad M_i \ge 0,$$ where $\mu$ is a certain
probability measure with unbounded support. For these polynomials,
we obtain the relative asymptotics with respect to orthogonal
polynomials related to $\mu$, Mehler--Heine type asymptotics and
their consequences about the asymptotic behaviour of the zeros.

To establish these results we use a new approach different from
the methods used in the literature up to now. The development of
this technique is highly motivated by the fact that the methods
used when $\mu$ is bounded do not work.

\end{abstract}

2000MSC: 42C05, 33C45.


\newpage

Corresponding author: Ana Pe\~{n}a, e-mail: anap@unizar.es

Departamento de Matem\'aticas.
Facultad de Ciencias.

C/ Pedro Cerbuna 12. Universidad de Zaragoza.

50009-Zaragoza  (Spain)

Fax: (+34) 976761338; Phone: (+34) 976761328

\newpage

\section{Introduction}

Let $\{\mu_i\}_{i=0}^r$  be Borel positive measures supported on the real line.
We can define the
Sobolev space:
\begin{align*}&W^{2,r}(\mu_0, \mu_1, \ldots, \mu_r):=\{f:
\int |f|^2d\mu_0 + \sum_{i=1}^r
\int|f^{(i)}|^2d\mu_i<+\infty \}
\end{align*} with the inner product
$$(f,g)=\int f\,g d\mu_0+\sum_{i=1}^r\int f^{(i)}\,g^{(i)} d\mu_i . $$

It is very well known that this inner product is nonstandard, that
is, $(xf,g)\ne (f, xg).$ Therefore, the nice properties of the
standard orthogonal polynomials, such as the three--term
recurrence relation, the interlacing properties of the zeros, etc,
do not hold any more. Then,  the powerful methods and techniques
developed for over a century to study \textit{standard} orthogonal
polynomials could not  work (in fact, they cannot work) for
Sobolev orthogonal polynomials. Then, a question arises:
\textit{is it necessary or interesting to study these
``pathological'' polynomials?} In our opinion the answer is
affirmative. It does not exist a general theory for these families
of orthogonal polynomials, for example,  up to now powerful tools
such as Rieman--Hilbert approach
 to obtain asymptotic properties of these
polynomials have not worked. This should be a motivation to pay
attention to Sobolev orthogonal polynomials, that is, to
investigate how to construct a more general theory in the same
sense as it was made for the standard orthogonal polynomials many
years ago. Furthermore, some applications of the Sobolev
orthogonality in the theory of standard orthogonal polynomials are
known, for instance, standard polynomials with nonstandard
parameters are not orthogonal in the usual sense but they are
orthogonal with respect to Sobolev inner products (see among
others \cite{aar}, \cite{cs} or \cite{kl}).

Thus, according to the above reasoning, in this paper we take a
step to get a better knowledge of the properties of the  Sobolev
orthogonal polynomials, more concretely of the discrete Sobolev
orthogonal polynomials. Let $\mu$ be  a finite positive Borel
measure supported on the real line, $c \in \mathbb{R}$ and $M_i\ge
0$ for  $i=0,1,...,r$. We consider an inner product of the form
\begin{equation*}
(f,g)=\int f(x)g(x)d\mu (x) +\sum_{i=0}^{r}M_i
f^{(i)}(c)g^{(i)}(c),
\end{equation*}
 and let $\{Q_{n}\}_{n\ge0}$ be  the corresponding sequence of monic orthogonal polynomials.
  More general products where cross--product terms appear in the discrete part have also been studied.
   But, recently  in \cite{kly} the authors prove that every symmetric bilinear form
  can be reduced to a diagonal case, that is, without cross--product terms.

 The aim is to compare the Sobolev orthogonal polynomials with the standard orthogonal
 polynomials associated with the measure $\mu$ in order to investigate how the addition of the
 derivatives in the inner product influences the orthogonal system.

 There exist many formal results for these polynomials: recurrence relation, location of zeros,
  differential formulas, and so on. However, little is known concerning the asymptotic properties. We want to
  remark that most of the general
  results have been obtained  when $\textrm{supp}(\mu)$ is a bounded set.  More precisely, in \cite{lmva}, the authors assume that $\mu$ is a measure
   for which the asymptotic behaviour
   of the orthogonal polynomials is known; the most relevant class of this type is the Nevai class $M(0,1)$ of orthogonal
    polynomials with appropriately converging recurrence  coefficients. They studied  the relative asymptotics when
     the mass point $c$ is outside the support of the measure.
     The same product  with the mass point in $\textrm{supp}(\mu)$ has been studied in \cite{rms}.

 In both papers, the key is the possibility to transform the
 Sobolev orthogonality into the standard quasi--orthogonality. As a consequence,
 we can express the polynomial $Q_n$ as a linear combination (with a fixed number of terms)
 of standard orthogonal polynomials $R_n$ corresponding  to the modified measure $d\nu=(x-c)^{r+1}d\mu$, that
 is,
 \begin{equation}\label{qo}
Q_{n}(x) = \sum_{j=0}^{r+1} a_{n}^j R_{n-j}(x) \,.
\end{equation}

We want to point out that in the bounded case, a straightforward
argument yields to prove  that all the connexion coefficients
$a_{n}^j$ are bounded. This behaviour of the coefficients together
with the fact that the orthogonal polynomials $R_n$ have an
adequate finite ratio asymptotics is enough reason
 to study  each term of (\ref{qo}) separately, in order to get the relative asymptotics
 for $Q_n$ (see \cite{lmva} and \cite{rms} where this technique is developed).

However, the situation is quite different if we deal with the
unbounded case. More concretely,  we consider the  Laguerre
probability  measure, that is,
$d\mu(x)=\frac{x^{\alpha}e^{-x}}{\Gamma(\alpha+1)}dx\,$ with
$\alpha>-1,$ and the inner product
\begin{equation} \label{intr-pr}
(f,g)_r=\frac{1}{\Gamma(\alpha+1)}\int_0^{\infty}
f(x)g(x)x^{\alpha}e^{-x}dx +\sum_{i=0}^{r}M_i
f^{(i)}(0)g^{(i)}(0),
\end{equation}
where $ M_i> 0, \, i=0,\dots,r $. Then we will see in Theorem \ref{corolario1} that the connexion
coefficients  which appear in (\ref{qo}) are unbounded. So, as we
will see later, when we try to obtain the relative asymptotics
with the techniques used for the bounded case and we take into
account the ratio asymptotics for Laguerre polynomials, we
come across a serious problem. Indeed, we find that each term of
(\ref{qo}) tends to infinity, all of them being the same order,
but with an alternating sign. Then, the idea that each term has a
limit does not work now, and therefore we can say that there is an
intrinsic difficulty related to the unbounded case.

 The interest lies in knowing the differences  between  the Laguerre polynomials and the
Sobolev  polynomials $Q_{n,r}$ orthogonal with respect to
(\ref{intr-pr}). Intuitively one can imagine that the asymptotic
differences in the complex plane should be around the perturbation
of the standard inner product involved in the Sobolev inner
product, that is, around the origin. To confirm this, first  we
get the relative asymptotics in Theorem \ref{relative-asymptotic}
and we prove that both families of orthogonal polynomials
$Q_{n,r}$ and $L_n^{\alpha}$ are identical asymptotically on
compact subsets of $\mathbb{C}\setminus [0,\infty).$

Later, we consider Mehler--Heine type formulas because they are
nice tools to describe the Laguerre--Sobolev type polynomials
around the origin. In \cite{alv-mb}, with $r=1$ and $M_0, M_1>0$
the authors find a behaviour pattern and they establish a
conjecture which  is reformulated in \cite{mm} more clearly than
in \cite{alv-mb}. That is:

\textit{If $M_i>0$ for $i=0, \ldots, r,$ in the inner product
(\ref{intr-pr}), then
\begin{equation*}
\lim_{n\to \infty} \frac{(-1)^n}{n!
n^{\alpha}}Q_{n,r}\left(\frac{x}{n+j}\right)=(-1)^{r+1}x^{-\alpha/2}J_{\alpha+2r+2}(2\sqrt{x}),
\end{equation*}
uniformly on compact subsets of the complex plane and on $j\in
\mathbb{N}\cup \{0\}$ where $\{Q_{n,r}\}_{n\ge0}$ is the sequence
of monic  orthogonal polynomials with respect to (\ref{intr-pr})
and $J_{\alpha}$ is the Bessel function of the first kind of order
$\alpha$.}

\medskip
In Theorem \ref{teorMH} we prove that this conjecture is true and
this is one of our main results. We would like to note that the
techniques used to prove it are not a simple generalization of the
ones used in \cite{alv-mb}. There, the authors consider the
algebraic expression

\begin{equation*}
Q_{n,1}(x)=B_0(n)L_n^{\alpha}(x)+B_1(n)xL_{n-1}^{\alpha+2}(x)+B_2(n)x^2L_{n-2}^{\alpha+4}(x)
\end{equation*}
where the coefficients $B_i(n)$ were given explicitly  in
\cite{km}.
 Now, in a discrete Laguerre--Sobolev inner product with an arbitrary number of terms, the problem
  is that we only have an algebraic expression given   in  \cite{k1}, but not
   the explicit expression of the coefficients $B_i(n)$, of which we only know that they  are a
    non trivial solution of a system with $r+1$ equations and $r+2$ unknowns.

The approach used to prove Theorem \ref{teorMH} is totally new. We
do not use algebraic tools but analytic ones. For this purpose, we
obtain a new and nice formula for the derivatives of $Q_{n,r}$ and
as a consequence, we achieve a uniform bound for the ratios
$\frac{Q_{n,r}^{(k)}(0)}{(L_n^{\alpha})^{(k)}(0)}$.

These Mehler--Heine type formulas are interesting twofold: they
provide the scaled asymptotics for $Q_{n,r}$ on compact sets of
the complex plane and they supply us with asymptotic information
about the location of the zeros of these polynomials in terms of
the zeros of other known special functions. More precisely,
applying Hurwitz's Theorem in a straightforward way, we prove that
there exists an acceleration of the convergence of $r+1$ zeros of
these Sobolev polynomials towards the origin.

Along the paper, we also deal with the possibility of some
$M_i=0$. We call  such Sobolev inner product, a
 Sobolev  inner product
with \textit{holes}. More concretely, we consider the inner
product
\begin{align*}
&(f,g)_{r,s} \\
&=\frac{1}{\Gamma(\alpha+1)}\int_{0}^{\infty}f(x)g(x)
x^{\alpha}e^{-x}dx+\sum_{i=0}^{r}M_i f^{(i)}(0)g^{(i)}(0)+M_s
f^{(s)}(0)g^{(s)}(0),
\end{align*}
where $s\ge r+2$ and $M_i>0$ for $i=0, \ldots, r$ and $i=s.$

 For this situation, we  also establish the relative asymptotics and the Mehler-Heine
type formulas for these  orthogonal polynomials. We want to remark
that this case has qualitative differences with respect to the
case without holes. For example, concerning the convergence
acceleration to $0$ of the zeros of the polynomials, the result
does not depend on the number of terms in the discrete part, but
it depends on the position of the hole. So,  despite  the presence
of the mass $M_s$, there only exists an acceleration of the
convergence of $r+1$ zeros such as it occurs in the case of the
inner product without holes.

As a consequence of all the previous results, using the
symmetrization process in the framework of Sobolev type orthogonal
polynomials, we can prove in Proposition \ref{hermite}  the
relative asymptotics and the Mehler--Heine type formulas  for the
generalized Hermite-Sobolev type polynomials. Furthermore, we hope
this method can be used with other unbounded measures for which we
have quite less explicit information about the corresponding
orthogonal polynomials.

The structure of the paper is as follows. In Section 2 we
introduce the notation, the basic tools, and some properties of
classical Laguerre polynomials. In Section 3 we obtain a new
auxiliary result in Laguerre--Sobolev type orthogonal polynomials
that we use to establish our main results in Sections 4 and 5.
Concretely, Section 4 is devoted to obtain the relative
asymptotics and the Mehler--Heine type formula for the orthogonal
polynomials with respect to a discrete Sobolev inner product with
positive masses and in Section 5 we get the corresponding ones for
the orthogonal polynomials with respect to an inner product  with
holes. We remark that exterior strong asymptotics for the
sequences of Sobolev orthogonal polynomials considered in Sections
4 and 5 are trivially deduced from the relative asymptotics
obtained here, because we know explicitly the exterior strong
asymptotics for the corresponding sequences of orthogonal
polynomials with respect to  $\mu.$ Furthermore, exterior
Plancherel--Rotach asymptotics for these Sobolev orthogonal
polynomials can also be obtained in a similar way with the
technique introduced here.

\bigskip

\section{Notation and basic results}

Along this work we will deal with classical Laguerre polynomials,
that is, polynomials orthogonal with respect
 to the inner product in the space of all polynomials with real coefficients
\begin{equation*}\label{innerproduct}
(p,q) = \frac{1}{\Gamma(\alpha +1)}\int_0^{\infty } p(x) q(x)\,
x^{\alpha} e^{-x} \,dx \,, \quad \alpha  > -1.
\end{equation*}
We will denote by $L_n^{\alpha}$ the $n$th monic Laguerre
polynomial.

Many of the properties of Laguerre polynomials can be seen, for
example, in  Szeg\H{o}'s book \cite{sz}. In what follows  we
summarize  those properties which will be used in this paper.

It is known that the \emph{monic Laguerre polynomials} are defined by
\begin{equation*}
L_n^{\alpha}(x) = (-1)^n n! \sum_{k=0}^n \left( \begin{matrix} n+\alpha \\ n-k  \end{matrix} \right)(-1)^k \frac{x^k}{k!}\,,
\end{equation*}
and their $L_2$-norm is
\begin{equation}\label{normaLn}
\Vert L_n^{\alpha} \Vert^2= \frac{1}{\Gamma(\alpha +1)}
\int_0^{\infty}
(L_n^{\alpha}(x))^2\,x^{\alpha}\,e^{-x}\,dx=\frac{\Gamma(n+\alpha+1)}{\Gamma(\alpha +1)}n! \,.
\end{equation}

The evaluation at $x=0$ of the polynomial $L_n^{\alpha}$ and its
successive derivatives are given by

\begin{equation}\label{derivadaLn}
(L_n^{\alpha})^{(k)}(0) = \frac{(-1)^k \, n!}{(n-k)!}
 \frac{\Gamma(\alpha +1)}{\Gamma(\alpha+k+1)} L_n^{\alpha}(0) =  \frac{(-1)^{n+k} \,
 n!}{(n-k)!}
 \frac{\Gamma(n+\alpha +1)}{\Gamma(\alpha+k+1)}.
\end{equation}

A useful tool to some estimates is the \textit{Stirling formula}:

\begin{equation*}
\Gamma(x+1) \sim x^x e^{-x} \sqrt{2\pi x} \quad ( x \to +\infty)
\end{equation*}
where the symbol $f(x) \sim g(x) \quad ( x \to a )$ stands for $\lim_{x \to a} \frac{f(x)}{g(x)} = 1.$
In particular,
\begin{equation}\label{equivgamma}
\Gamma(n+\alpha+1) \sim n! \, n^{\alpha}.
\end{equation}
As a consequence, from (\ref{normaLn}), (\ref{derivadaLn}), and
(\ref{equivgamma}), we get
\begin{equation}\label{estimate}
\lim_n \frac{(L_n^{\alpha}(0))^2}{\Vert L_n^{\alpha}\Vert^2 \, n^{\alpha}} =
\frac{1}{\Gamma(\alpha +1)} \,.
\end{equation}

The following asymptotic results are known. They can be deduced
from Perron's formula in Szeg\H{o}'s book \cite{sz},
\begin{equation}\label{asintrelatL}
\lim_n \frac{n L_{n-1}^{\alpha}(x)}{L_n^{\alpha}(x)} = -1,
\end{equation}

\begin{equation}\label{asintrelatLdistintosparametros}
\lim_n \frac{n^{1/2} L_{n}^{\alpha}(x)}{L_n^{\alpha +1}(x)} = \sqrt{-x},
\end{equation}
both uniformly on compact subsets of $\mathbb{C} \setminus
[0,\infty)$.

\bigskip

The $n$th kernel for the Laguerre polynomials
$K_n(x,y)=\displaystyle \sum _{i=0}^n \frac{L_i^{\alpha}(x)L_i^{\alpha}(y)}{\Vert
L_i^{\alpha} \Vert ^2}$ satisfies the Christoffel--Darboux formula
\begin{equation*}
K_{n}(x,y)=\frac{1}{{\Vert L_n^{\alpha}
\Vert}^2}\frac{L_{n+1}^{\alpha}(x)L_n^{\alpha}(y)-L_{n+1}^{\alpha}(y)L_n^{\alpha}(x)}{x-y} \,.
\end{equation*}
As usual, we denote the derivatives of the  kernels  by
$$K_n^{(k,s)}(x,y)=\displaystyle\frac{\partial^{k+s}}{\partial{x^k}
\partial{y^s}}K_n(x,y)
=\displaystyle \sum _{i=0}^n
\frac{(L_i^{\alpha})^{(k)}(x)(L_i^{\alpha})^{(s)}(y)}{\Vert L_i^{\alpha} \Vert ^2}$$
with $k,s  \in \mathbb{N}\cup \{0\}$ and the convention $K_n^{(0,0)}(x,y)=K_n(x,y)$.

\bigskip

In the next lemma we show some formulas for the derivatives of the
kernels that we will need throughout the paper.

\begin{lemma} \label{nucleos}
The derivatives of the kernels of the Laguerre polynomials, for
$k,s  \in \mathbb{N}\cup \{0\}$, satisfy
\begin{itemize}
\item[(a)]
\begin{equation*}
K_{n-1}^{(0,s)}(x,0)=\frac{1}{{\Vert L_{n-1}^{\alpha} \Vert}^2}\frac{s!}{x^{s+1}}\left[
P_s(x,0;L_{n-1}^{\alpha})L_{n}^{\alpha}(x)-P_s(x,0;L_{n}^{\alpha})L_{n-1}^{\alpha}(x) \right]
\end{equation*} where
$P_s(x,0;f)$ is the $\it{s}$th Taylor polynomial of $f$ at $0$.
\end{itemize}

\begin{itemize}
\item[(b)]

\begin{equation*}
K_{n-1}^{(k,s)}(0,0)=\frac{k!\, s!}{\Vert L_{n-1}^{\alpha}\Vert^2}
\sum_{j=0}^s \frac{k+s+1-2j}{n-j}\, \frac{(L_{n-1}^{\alpha})^{(j)}(0)(L_{n}^{\alpha})^{(k+s+1-j)}(0)}{j! \, (k+s+1-j)!}\, ,
\end{equation*}

\begin{equation*}
K_{n-1}^{(k,0)}(0,0)=(-1)^k \frac{\Gamma(\alpha +n+1)}{(n-(k+1))!\, \Gamma(\alpha+k+2)}\,.
\end{equation*}

\end{itemize}
\end{lemma}
\textbf{Proof.} (a) The result follows from the Christoffel-Darboux formula
and Leibniz's rule.

(b) Observe that, according to Taylor's formula,
$\frac{1}{k!}K_{n-1}^{(k,s)}(0,0)$ is precisely  the coefficient of
$x^k$  in $K_{n-1}^{(0,s)}(x,0)$, therefore

\begin{align*}
&K_{n-1}^{(k,s)}(0,0)\\
&=\frac{k! \, s!}{\Vert L_{n-1}^{\alpha}\Vert^2}
\sum_{j=0}^s \frac{(L_{n-1}^{\alpha})^{(j)}(0)(L_{n}^{\alpha})^{(k+s+1-j)}(0)-(L_{n}^{\alpha})^{(j)}(0)(L_{n-1}^{\alpha})^{(k+s+1-j)}(0)}{j! \, (k+s+1-j)!}.
\end{align*}
In particular for $s=0$, straightforward computations lead us to
conclude this lemma. $\quad \Box$

\medskip

Along the paper we work with sequences of monic orthogonal
polynomials, and we use the acronym SMOP for them.

\section {Auxiliary results}

From now on $\{ Q_{n,r} \}_{n\ge0}$ denotes the sequence of monic
Laguerre--Sobolev orthogonal polynomials  with respect to an inner
product of the form
\begin{equation}\label{pr-ls1}
(p,q )_r = \frac{1}{\Gamma(\alpha +1)}\int_0^{\infty} p(x)q(x) \,
x^{\alpha} e^{-x} \,dx + \sum_{i=0}^r M_i p^{(i)}(0) q^{(i)}(0),
\end{equation}
where $\alpha  > -1$ and  $\, M_i > 0, \, i = 0, \dots, r.$
Notice that all the masses in the discrete part of this inner
product are positive.

We write $K_{n,r}$ for the corresponding $n$th kernel, that is
$\displaystyle K_{n,r}(x,y)= \sum _{j=0}^n
\frac{Q_{j,r}(x)Q_{j,r}(y)}{(Q_{j,r},Q_{j,r})_r} \,,$ and
$K_{n,r}^{(k,s)}$ for the derivatives of the kernels.

Observe that, in fact, $(.,. )_r$, $Q_{n,r} $, $K_{n,r}$ and
$K_{n,r}^{(k,s)}$ also depend  on the parameter $\alpha$ but for
simplicity we have omitted it in the notations.

In the next lemma, we obtain an asymptotic estimation for
$Q_{n,r}^{(k)}(0)$, $k\ge 0$, that will play  an important role
along this paper. To do this, we need to know the \lq\lq size" of
the kernels of $Q_{n,r}$ and their derivatives. The following
discrete version of l'Hospital's rule (see e.g. \cite{k}) will be
very helpful to easily calculate some limits:

\textit{Stolz Criterion.} Let $\{x_n\}$ and $\{y_n\}$ be real sequences. Suppose that $\{y_n\}$ is monotonic and $y_n \not= 0$ for all $n$. If
$\lim_n (x_{n+1} - x_n)/(y_{n+1} - y_n) = L \in \mathbb{R} \cup \{\pm \infty\}$ exists, then $\lim_n x_n/y_n = L$ provided either
$\lim_n x_n = \lim_n y_n = 0$ or $\lim_n y_n = \pm \infty$.

\begin{lemma}\label{c}
Let $ Q_{n,r} $ be the monic polynomials orthogonal with respect to the inner product (\ref{pr-ls1}). Then the following statements hold:
\begin{itemize}
\item[(a)] For $0 \le k \le r$,
$$ \frac{Q_{n,r}^{(k)}(0)}{(L_n^{\alpha})^{(k)}(0)} \sim \frac{C_{r,k}}{n^{\alpha+2k+1}} ,$$

where $C_{r,k}$ is a nonzero real number independent of $n$.

For $k \ge r+1$,
$$\lim_n \frac{Q_{n,r}^{(k)}(0)}{(L_n^{\alpha})^{(k)}(0)}=
\frac{k!}{(k-(r+1))!}\frac{\Gamma(\alpha+k+1)}{\Gamma(\alpha+r+k+2)}\,  .$$

\item[(b)]
$$\lim_n \frac{(Q_{n,r},Q_{n,r})_r}{\Vert L_n^{\alpha} \Vert ^2} = 1 \,.$$
\end{itemize}
\end{lemma}
\textbf{Proof.}  We use mathematical induction on $r \in
\mathbb{N} \cup \{0\}$.

If $r=0$, the Fourier expansion of the polynomial $Q_{n,0}$ in the orthogonal basis $\{L_n^{\alpha}\}_{n\ge0}$ leads to
\begin{equation*}
Q_{n,0}(x) = L_n^{\alpha}(x) - M_0 Q_{n,0}(0) K_{n-1}(x,0) \,,
\end{equation*}
and therefore
\begin{equation}\label{Qn0(x)}
Q_{n,0}(x) = L_n^{\alpha}(x) - \frac{M_0 L_n^{\alpha}(0)}{1+ M_0 K_{n-1}(0,0)} K_{n-1}(x,0) \,.
\end{equation}
As a consequence of (\ref{derivadaLn}) and Lemma \ref{nucleos}
(b), we obtain $(a)$ for $r=0$.

\noindent Using (\ref{Qn0(x)}), we have

\begin{equation*}
(Q_{n,0},Q_{n,0})_0 = \Vert L_n^{\alpha} \Vert ^2 +\frac{M_0 (L_n^{\alpha}(0))^2}{1+M_0 K_{n-1}(0,0)}\,.
\end{equation*}
Thus, from (\ref{estimate}) and  Lemma \ref{nucleos} (b) it
follows $(b)$ for $r=0$.

Suppose now that $(a)$ and $(b)$ hold for the SMOP $\{ Q_{n,r}
\}_{n\ge0}$ with $r>0$, then we are going to deduce that they  are
also true for the sequence $\{Q_{n,r+1}\}_{n\ge0}$. To do this, we
observe that
\begin{equation}\label{prodescalar}
(p,q)_{r+1} = (p,q)_r + M_{r+1} p^{(r+1)}(0) q^{(r+1)}(0)
\end{equation}
and therefore
\begin{equation}\label{formulasencilla}
Q_{n,r+1}(x) =Q_{n,r}(x) - M_{r+1}Q_{n,r+1}^{(r+1)}(0)K_{n-1,r}^{(0,r+1)}(x,0)  \,.
\end{equation}
Taking derivatives $r+1$ times in (\ref{formulasencilla})  and
evaluating at $x=0$, we obtain
\begin{equation}\label{Qnr(r+1)}
Q_{n,r+1}^{(r+1)}(0) =\frac{Q_{n,r}^{(r+1)}(0)}{1 + M_{r+1}K_{n-1,r}^{(r+1,r+1)}(0,0)} \,.
\end{equation}

\noindent Taking now derivatives  $k$ times in
(\ref{formulasencilla}), evaluating at $x=0, $ and using
(\ref{Qnr(r+1)}) we get
\begin{equation}\label{Qnr/Ln}
\frac{Q_{n,r+1}^{(k)}(0)}{(L_n^{\alpha})^{(k)}(0)} =\frac{Q_{n,r}^{(k)}(0)}{(L_n^{\alpha})^{(k)}(0)} - \frac{M_{r+1}K_{n-1,r}^{(k,r+1)}(0,0)}{1+M_{r+1}K_{n-1,r}^{(r+1,r+1)}(0,0)}
\frac{Q_{n,r}^{(r+1)}(0)}{(L_n^{\alpha})^{(k)}(0)} \,.
\end{equation}

\medskip
Before taking limits in the last expression, we need to estimate $K_{n-1,r}^{(k,r+1)}(0,0)$. Applying Stolz criterion,
 the induction hypothesis for $\{ Q_{n,r} \}_{n\ge0}$, (\ref{derivadaLn}) and (\ref{estimate}), we obtain for $k\ge r+1$

\begin{align}\label{limKsobren}
\nonumber
& \lim_n \frac{K_{n-1,r}^{(k,r+1)}(0,0)}{n^{\alpha+k+r+2}}
=\lim_n  \frac{Q_{n-1,r}^{(k)}(0)Q_{n-1,r}^{(r+1)}(0)}{\Vert L_{n-1}^{\alpha}  \Vert^2(\alpha+k+r+2) \,\, n^{\alpha+k+r+1}} \\ \nonumber
&=\frac{(-1)^{k+r+1}\,\Gamma(\alpha+1)}{(\alpha+k+r+2) \Gamma(\alpha+k+1) \Gamma(\alpha+r+2)} \lim_n \left[ \frac{Q_{n-1,r}^{(k)}(0)}{(L_{n-1}^{\alpha})^{(k)}(0)} \,\, \frac{Q_{n-1,r}^{(r+1)}(0)}{(L_{n-1}^{\alpha})^{(r+1)}(0)} \right] \\
&=\frac{k! (r+1)!}{(k-(r+1))!}\frac{(-1)^{k+r+1}\Gamma(\alpha+1)}{\Gamma(\alpha+k+r+3)\Gamma(\alpha+2r+3)} \,,
\end{align}

\noindent and therefore, from (\ref{Qnr/Ln}), we get $(a)$ for $k \ge r+1$.

Now, if $0 \le k \le r$, to estimate the size of
$K_{n-1,r}^{(k,r+1)}(0,0)$, we use  Stolz criterion again and
thus, we obtain

\begin{equation*}
\lim_n \frac{K_{n-1,r}^{(k,r+1)}(0,0)}{n^{r+1-k}}= (-1)^{r+1+k} \frac{(r+1)!}{r+1-k} \frac{C_{r,k}\, \Gamma(\alpha+1)}{\Gamma(\alpha+k+1)\Gamma(\alpha+2r+3)} \,.
\end{equation*}

Using the induction hypothesis and substituting all these results
in the right--hand side of (\ref{Qnr/Ln}) we get
$$ \frac{Q_{n,r+1}^{(k)}(0)}{(L_n^{\alpha})^{(k)}(0)} \sim \frac{C_{r+1, k}}{n^{\alpha+2k+1}}  \, ,$$
where $C_{r+1, k} = - \frac{\alpha+k+r+2}{r+1-k} C_{r,k} \not= 0$. Therefore the proof of $(a)$ is complete.

To finish the proof of the Lemma, we only need to prove $(b)$ for
$\{ Q_{n,r+1} \}_{n\ge0}$. As in the case $r=0$, from
(\ref{prodescalar}) and (\ref{Qnr(r+1)}),  we get
\begin{equation*}
(Q_{n,r+1},Q_{n,r+1})_{r+1} =
(Q_{n,r},Q_{n,r})_r + \frac{M_{r+1} (Q_{n,r}^{(r+1)}(0))^2}{1+ M_{r+1} K_{n-1,r}^{(r+1,r+1)}(0,0)} \,.
\end{equation*}
 Using (a) for $k=r+1$, (\ref{normaLn}), (\ref{derivadaLn}) and (\ref{limKsobren}) we achieve the result.
$ \quad \Box$

\section{Main results}

 As we have mentioned in the introduction, if we consider a general discrete Sobolev inner product where
  the support of the measure $\mu$ is a bounded set, the key
  used to obtain some results is the possibility to transform  the
  Sobolev orthogonality into a standard  quasi--orthogonality.

 Now, in our particular case, the sequence $\{ Q_{n,r}\}_{n\ge0}$  orthogonal with
respect to the inner product defined by (\ref{pr-ls1}) is
quasi--orthogonal of order $r+1$ with respect to the Laguerre
weight $x^{\alpha+r+1} e^{-x},$ that is, $$\int_0^{+\infty} p(x)
Q_{n,r}(x) x^{\alpha+r+1}e^{-x}dx=0,$$ for every polynomial $p$
with $\deg p\le n-(r+1)-1$. Therefore, as an immediate consequence
we have a \emph{connexion formula} of the form
\begin{equation}\label{conexion}
Q_{n,r}(x) = \sum_{j=0}^{r+1} a_{n,r}^j L_{n-j}^{\alpha+r+1}(x)
\,, \quad a_{n,r}^0=1.
\end{equation}

 In this section, our first
 effort is devoted to obtain   the size of the \emph{connexion coefficients}  $a_{n,r}^j$. We introduce a fruitful
   and new technique which  leads us to know the asymptotic behaviour  of these coefficients.

Then, we get two asymptotic properties of the Sobolev type
orthogonal polynomials $Q_{n,r}$. The first one is the
\emph{relative asymptotics}, which shows that the
Laguerre--Sobolev type polynomials are very  similar to the
Laguerre polynomials outside the support of the measure. The
second one is the so--called Mehler--Heine type formula which
shows  how the presence of the masses
  in  the inner product changes the asymptotic behaviour around  the
  origin.

As we will see later, it is worth noticing that the knowledge of
the asymptotic behaviour of the connexion coefficients is not
enough  to obtain both asymptotic properties. Thus, we can assure
that the  techniques used in the bounded case  do not work now.

\subsection{Connexion coefficients}

 Notice that the sequence $\{ Q_{n,r}\}_{n\ge0}$ is quasi--orthogonal of order $r+1$ with respect to the Laguerre weight $x^{\alpha+r+1} e^{-x}$ and therefore we have a \emph{connexion formula} of the form
\begin{equation*}
Q_{n,r}(x) = \sum_{j=0}^{r+1} a_{n,r}^j L_{n-j}^{\alpha+r+1}(x) \,.
\end{equation*}

Using this expansion, in the next lemma we obtain a new algebraic
expression with a
 nice  structure which relates the derivative of order $k+1$ of the polynomials $Q_{n,r}$ to
  the derivative of order $k$.

\begin{lemma} \label{lema2}
Fixed $r \ge 1$, let $\{ Q_{n,r}\}_{n \ge 0}$ be the SMOP with respect to the inner product (\ref{pr-ls1}).  Then, we have for $0\le k \le n-1$,
\begin{align}\label{ratioderiv}
& \frac{Q_{n,r}^{(k+1)}(0)}{(L_n^{\alpha})^{(k+1)}(0)}  =
\frac{\alpha+k+1}{\alpha+r+k+2}\frac{Q_{n,r}^{(k)}(0)}{(L_n^{\alpha})^{(k)}(0)}   \\ \nonumber
& +   \frac{\Gamma(\alpha+r+2)}{\Gamma(\alpha+r+k+3)} \sum_{i=1}^{k+1} \left( \begin{matrix} k \\ i-1
 \end{matrix} \right) \frac{\Gamma(\alpha+k+2)}{\Gamma(\alpha+i+1)} \frac{A_{n,r}^i}{(L_n^{\alpha})^{(i)}(0)}
 \,,
\end{align}

\noindent where
\begin{equation*}
A_{n,r}^i= \sum_{j=i}^{r+1} \frac{j!}{(j-i)!}\,a_{n,r}^j \, L_{n-j}^{\alpha+r+1}(0), \quad i=0, 1, \dots, r+1 \,,
\end{equation*}
and the coefficients $a_{n,r}^j$ are  those of (\ref{conexion}).
By convention,  $A_{n,r}^i =0 ,$ when $i > r+1$. Besides,
\begin{equation}\label{limcorchete}
\lim_n \frac{A_{n,r}^i}{(L_n^{\alpha})^{(i)}(0)} = \begin{cases}  0 & \text {if $\quad 0 \le i \le r$} \, ,\\
(r+1)! & \text {if $\quad i = r+1$} \,.\end{cases}
\end{equation}

\end{lemma}
\textbf{Proof.} Taking derivatives  $k+1$ times in
(\ref{conexion}), evaluating at $x=0$, and
 using  (\ref{derivadaLn}) several times we get, for $k \ge 0$,

\begin{align*}
&Q_{n,r}^{(k+1)}(0) = \sum_{j=0}^{r+1} a_{n,r}^j (L_{n-j}^{\alpha+r+1})^{(k+1)}(0) = - \sum_{j=0}^{r+1} a_{n,r}^j \frac{n-j-k}{\alpha+r+k+2} (L_{n-j}^{\alpha+r+1})^{(k)}(0) \\
& = - \frac{n-k}{\alpha+r+k+2} Q_{n,r}^{(k)}(0) + \frac{1}{\alpha+r+k+2} \sum_{j=1}^{r+1} j a_{n,r}^j (L_{n-j}^{\alpha+r+1})^{(k)}(0) \\
& = - \frac{n-k}{\alpha+r+k+2} Q_{n,r}^{(k)}(0) \\
&+ \frac{(-1)^{k}\,\Gamma(\alpha+r+2)}{\Gamma(\alpha+r+k+3)} \sum_{j=1}^{r+1} j a_{n,r}^j \frac{(n-j)!}{(n-j-k)!} L_{n-j}^{\alpha+r+1}(0)\,.
\end{align*}
According to formula (5) in page 8 of \cite{r},
\begin{equation*}
\frac{(n-j)!}{(n-j-k)!} = k! \sum_{i=0}^k (-1)^{i} \left( \begin{matrix} j-1 \\ i  \end{matrix} \right)
 \left( \begin{matrix} n-1-i \\ k-i  \end{matrix} \right),
\end{equation*}
and therefore
\begin{equation*}
\frac{(n-j)!}{(n-j-k)!} = \frac{(j-1)!}{(n-(k+1))!} \sum_{i=1}^{k+1} (-1)^{i-1} \left( \begin{matrix} k \\ i-1  \end{matrix} \right) \frac{(n-i)!}{(j-i)!}\,.
\end{equation*}
Thus, the above expression can be written in the form:

\begin{align*}
& Q_{n,r}^{(k+1)}(0) =
- \frac{n-k}{\alpha+r+k+2} Q_{n,r}^{(k)}(0)  \\
& +  \frac{ (-1)^{k}\Gamma(\alpha+r+2)}{\Gamma(\alpha+r+k+3)} \sum_{i=1}^{k+1} (-1)^{i-1}\left( \begin{matrix} k \\ i-1  \end{matrix} \right) \frac{(n-i)!}{(n-(k+1))!} A_{n,r}^i\,,
\end{align*}

\noindent which leads to
\begin{align*}
&  \frac{Q_{n,r}^{(k+1)}(0)}{(L_n^{\alpha})^{(k+1)}(0)}  = - \frac{n-k}{\alpha+r+k+2} \frac{(L_n^{\alpha})^{(k)}(0)}{(L_n^{\alpha})^{(k+1)}(0)}
\frac{Q_{n,r}^{(k)}(0)}{(L_n^{\alpha})^{(k)}(0)} \\
&+ \frac{ (-1)^{k}\Gamma(\alpha+r+2)}{\Gamma(\alpha+r+k+3)} \sum_{i=1}^{k+1} (-1)^{i-1}\left( \begin{matrix} k \\ i-1  \end{matrix} \right) \frac{(n-i)!}{(n-(k+1))!} \frac{(L_n^{\alpha})^{(i)}(0)}{(L_{n}^{\alpha})^{(k+1)}(0)}\frac{A_{n,r}^i}{(L_n^{\alpha})^{(i)}(0)} \\
&=\frac{\alpha+k+1}{\alpha+r+k+2} \frac{Q_{n,r}^{(k)}(0)}{(L_n^{\alpha})^{(k)}(0)}    \\
& +   \frac{\Gamma(\alpha+r+2)}{\Gamma(\alpha+r+k+3)}  \sum_{i=1}^{k+1}  \left( \begin{matrix} k \\ i-1  \end{matrix} \right) \frac{\Gamma(\alpha+k+2)}{\Gamma(\alpha+i+1)} \frac{A_{n,r}^i}{(L_n^{\alpha})^{(i)}(0)}\,,
\end{align*}
where we have used  expression (\ref{derivadaLn}). So, the first part of the lemma is proved.

To prove (\ref{limcorchete}),  it is enough to apply  a recursive
procedure beginning with $k=0$ in (\ref{ratioderiv}) and  take
into account Lemma \ref{c} (a). $ \quad \Box$

\medskip

The above result allows us to obtain an estimate of the connexion
coefficients $a_{n,r}^j$. Observe that the following result is the
break point with respect to the techniques used in the bounded
case because we  establish that the coefficients in
(\ref{conexion}) are unbounded for $j\ge 1.$

\begin{theorem}\label{corolario1}
Let $a_{n,r}^j$ be the connexion coefficients which appear in
(\ref{conexion}). Then, we have
\begin{equation*}
\underset{n}\lim\frac{a_{n,r}^j }{n^j}=\left( \begin{matrix} r+1 \\ j  \end{matrix} \right)  \,, \quad 0 \le j \le r+1\,.
\end{equation*}
\end{theorem}
\textbf{Proof.} From (\ref{limcorchete}) and the expression
\begin{equation*}
 A_{n,r}^{r+1}= (r+1)! \, a_{n,r}^{r+1}  L_{n-r-1}^{\alpha+r+1}(0) = \frac{(r+1)! \,\, a_{n,r}^{r+1}}{n(n-1)\dots(n-r)} (L_n^{\alpha})^{(r+1)}(0)\,,
\end{equation*}
 it follows easily
$$\underset{n}\lim\frac{a_{n,r}^{r+1}}{ n^{r+1}}=1.$$

 A recurrence procedure leads to the result. Indeed, we assume that the result holds for $k+1\le j \le r+1$ and
 we will show that it is true for $j=k$.
 From (\ref{limcorchete}) for $i=k$ we can obtain

\begin{align*}
\underset{n}\lim  \sum_{j=k}^{r+1}  \frac{j!}{k!\,(j-k)!} \frac{a_{n,r}^j}{n^{r+1-k}} \frac{L_{n-j}^{\alpha+r+1}(0)}{(L_n^{\alpha})^{(k)}(0)}=0\,.
\end{align*}
From (\ref{derivadaLn}) and (\ref{equivgamma}), we get

\begin{align*}
\underset{n}\lim  \sum_{j=k}^{r+1} (-1)^{k-j} \left( \begin{matrix} j \\ k  \end{matrix} \right) \frac{a_{n,r}^j}{n^{j}} =0.
\end{align*}
From the assumption
\begin{equation*}
\underset{n}\lim\frac{a_{n,r}^j }{n^j}=\left( \begin{matrix} r+1 \\ j  \end{matrix} \right)  \,, \quad k+1 \le j \le r+1,\,
\end{equation*}
\noindent we have
\begin{align*}
\underset{n}\lim \frac{a_{n,r}^k}{n^{k}} &= \sum_{j=k+1}^{r+1} (-1)^{k+1-j} \left( \begin{matrix} j \\
k  \end{matrix} \right) \left( \begin{matrix} r+1 \\ j  \end{matrix} \right)\\ &= \left( \begin{matrix} r+1
 \\ k  \end{matrix} \right)\sum_{j=k+1}^{r+1} (-1)^{k+1-j} \left( \begin{matrix} r+1-k\\ j-k
 \end{matrix} \right)= \left( \begin{matrix} r+1 \\ k  \end{matrix}
 \right). \quad \Box
\end{align*}

\bigskip

\subsection{Relative asymptotics}

The relative asymptotics  for the Sobolev type orthogonal
polynomials with respect to the Laguerre polynomials cannot be
deduced as a consequence of the connexion formula
(\ref{conexion}). Indeed, from this formula we have

 $$\frac{Q_{n,r}(x)}{L_n^{\alpha}(x)}=\sum_{j=0}^{r+1} a_{n,r}^j \frac{L_{n-j}^{\alpha+r+1}(x)}{L_n^{\alpha}(x)}.$$

Applying  Theorem \ref{corolario1}, formulas (\ref{asintrelatL})
and (\ref{asintrelatLdistintosparametros}) in the above expression
we deduce that each term tends to infinity with the same order but
with an alternating sign, that is,
$$a_{n,r}^j \frac{L_{n-j}^{\alpha+r+1}(x)}{L_n^{\alpha}(x)}\sim (-1)^j\,\left( \begin{matrix} r+1 \\ j  \end{matrix} \right)\,\left(\frac{1}{\sqrt{-x}}\right)^{r+1}
\,n^{\frac{r+1}{2}},$$
\noindent uniformly on compact subsets of $\mathbb{C} \setminus [0,\infty)$.

\bigskip

It is important to remark that  the asymptotic behaviour  of the
above terms  is totally different from  the one when we consider
the  measures in the Nevai class (bounded case) in which every
term has a finite limit.

\begin{theorem}\label{relative-asymptotic}

Let $\{ Q_{n,r} \}_{n\ge0}$ be the SMOP with respect to the inner product defined by (\ref{pr-ls1}).
 Then
\begin{equation*}
\lim_n \frac{Q_{n,r}(x)}{L_n^{\alpha}(x)} = 1,
\end{equation*}
uniformly on compact subsets of $\mathbb{C} \setminus [0,\infty)$.
\end{theorem}
\textbf{Proof.} From the Fourier expansion of the polynomial
$Q_{n,r}$ in terms of the Laguerre polynomials we have
\begin{equation*}
 \frac{Q_{n,r}(x)}{L_n^{\alpha}(x)} = 1 - \sum_{i=0}^r M_i Q_{n,r}^{(i)}(0)  \frac{K_{n-1}^{(0,i)}(x,0)}{L_n^{\alpha}(x)}  \,.
\end{equation*}
 Substituting  Lemma \ref{nucleos} (a) in this  formula, we get
\begin{equation*}
 \frac{Q_{n,r}(x)}{L_n^{\alpha}(x)} = 1 - \sum_{i=0}^r \frac{M_i \, i!}{\Vert L_{n-1}^{\alpha}\Vert^2 \, x^{i+1}Q_{n,r}^{(i)}(0) P_i(x,0; L_{n-1}^{\alpha})}
 \left[1- \frac{P_i(x,0; L_n^{\alpha})}{P_i(x,0; L_{n-1}^{\alpha})} \frac{L_{n-1}^{\alpha}(x)}{L_n^{\alpha}(x)} \right] \,,
\end{equation*}
where $P_i(x,0;f)$ is the $i$-th Taylor polynomial of $f$ in $0$.

Now we analyze each one of the terms of this sum. Since
\begin{equation}\label{equivT}
\lim_n \frac{P_i(x,0; L_n^{\alpha})}{(L_n^{\alpha})^{(i)}(0)}= \frac{x^i}{i!} \,,
\end{equation}

\noindent from (\ref{derivadaLn}) and (\ref{asintrelatL}), we have, for $i=0,\dots,r$,
$$ \lim_n \left[1- \frac{P_i(x,0; L_n^{\alpha})}{P_i(x,0; L_{n-1}^{\alpha})}
 \frac{L_{n-1}^{\alpha}(x)}{L_n^{\alpha}(x)} \right] = 0,$$
uniformly on compact subsets of $\mathbb{C} \setminus [0,\infty)$.

Moreover, taking into account (\ref{derivadaLn}), (\ref{equivT}),
and Lemma \ref{c} (a), there exists
\begin{equation*}
\lim_n \frac{M_i \, i! Q_{n,r}^{(i)}(0) P_i(x,0; L_{n-1}^{\alpha})}{\Vert L_{n-1}^{\alpha}\Vert^2 \, x^{i+1}} \in \mathbb{C}\,.
\end{equation*}
Therefore each one of the terms in the sum tends to $0$ uniformly  on compact subsets of
 $\mathbb{C} \setminus [0,\infty)$ and the result follows.
$\quad \Box$

\subsection{Mehler--Heine type formulas}

Once we have proved in the previous subsection that both sequences
of orthogonal polynomials, $\{ Q_{n,r}\}_{n\ge 0}$ and
$\{L_n^{\alpha} \}_{n\ge 0},$ are asymptotically identical on
compact subsets of $\mathbb{C} \setminus [0,\infty),$  we
establish their differences through Mehler--Heine type formulas
which describe the asymptotic behaviour around the origin. First
of all, we recall the corresponding formula for the  monic
Laguerre polynomials, (see \cite[Th.8.1.3]{sz}):
\bigskip

\begin{equation}\label{MH Ln} \lim_n
\frac{(-1)^n}{n! \,n^{\alpha}}
L_n^{\alpha}\left(\frac{x}{n+j}\right)=x^{- \alpha/2} J_{\alpha}(2
\sqrt{x}),
\end{equation}
uniformly  on  compact  subsets of $\mathbb{C}$ and uniformly on $j\in \mathbb{N}\cup \{0\}$,
where $J_{\alpha}$  is the Bessel function of the first kind of order
$\alpha$ ($\alpha >-1$), defined by
\begin{equation*}
J_{\alpha}(x)=\sum_{n=0}^{\infty} \frac{(-1)^n}{n!
\,\Gamma(n+\alpha +1)} \left( \frac{x}{2} \right)^{2n+\alpha}.
\end{equation*}
Now, we want to obtain a similar result for the SMOP $\{ Q_{n,r}\}_{n\ge0}$.

As it occurs in the study  of the relative asymptotics, the
Mehler--Heine type formulas cannot be deduced as a consequence of
the connexion formula (\ref{conexion}). Indeed, from this formula
we have

\begin{equation*}
 \frac{(-1)^n}{n!\,n^{\alpha}} Q_{n,r}\left(\frac{x}{n+j}\right)=
 \sum_{i=0}^{r+1} a_{n,r}^i \frac{(-1)^n}{n!\,n^{\alpha}}
 L_{n-i}^{\alpha+r+1}\left(\frac{x}{n+j}\right).
\end{equation*}

Therefore, applying Theorem \ref{corolario1} and the Mehler--Heine
type formula for Laguerre polynomials, we deduce that each term
tends to infinity with the same order but with an alternating
sign.

Thus, to get the result for $\{ Q_{n,r}\}_{n\ge0}$, we focus on
the problem in a different way. We write the Taylor expansion of
the polynomial $ Q_{n,r}$

\begin{equation*}
 \frac{(-1)^n}{n!\,n^{\alpha}} Q_{n,r}\left(\frac{x}{n+j}\right) = \sum_{k=0}^n \frac{(-1)^n}{n!\,n^{\alpha}}  \frac{Q_{n,r}^{(k)}(0)}{(L_n^{\alpha})^{(k)}(0)}
 \frac{(L_n^{\alpha})^{(k)}(0)}{k!} \frac{x^k}{(n+j)^k}\,.
 \end{equation*}

Then, to calculate the limit, we are going to use  the Lebesgue's
dominated convergence theorem. For this purpose, we need to find a
uniform bound   for  the ratios
  $\displaystyle
  \frac{Q_{n,r}^{(k)}(0)}{(L_n^{\alpha})^{(k)}(0)}.$ It is clear
  that  when
    we take derivatives $k$ times in the expression (\ref{conexion})
    the connexion coefficients  do not change. Then, it could be thought
     about the possibility  to obtain this uniform bound from this formula. But again we come across
      the same problem, more concretely

\begin{equation*}
 \frac{Q_{n,r}^{(k)}(0)}{(L_n^{\alpha})^{(k)}(0)}=\sum_{i=0}^{r+1} a_{n,r}^i
  \frac{
  (L_{n-i}^{\alpha+r+1})^{(k)}(0)}{(L_n^{\alpha})^{(k)}(0)},
 \end{equation*}

\noindent and therefore, from Theorem \ref{corolario1}, again each
term tends to infinity with the same order, concretely $n^{r+1}$,
but with an alternating sign.

Then, we come back to the useful expression  established  in Lemma
\ref{lema2} which leads us to obtain a uniform bound for the
ratios
  $\displaystyle \frac{Q_{n,r}^{(k)}(0)}{(L_n^{\alpha})^{(k)}(0)}.$

\begin{lemma} \label{lema3}
Let $\{ Q_{n,r}\}_{n \ge 0}$ be the SMOP with respect to the inner product (\ref{pr-ls1}).
 Then, fixed $r \ge 1$ there exists a positive integer number $n_0$ such that for all $n \ge n_0$ and
 for all $k$ with $ r+1\le k \le n,$ the inequality
\begin{equation*}
\left| \frac{Q_{n,r}^{(k)}(0)}{(L_n^{\alpha})^{(k)}(0)} \right|  \le 2 (r+1) \frac{k!}{(k-r)!} (k-(r-1)) \frac{\Gamma(\alpha+k+1)}{\Gamma(\alpha+r+k+2)} \,,
\end{equation*}
 holds. Furthermore,  for $r\ge0$, there exists
$n_0 \in \mathbb{N}$ such that,
\begin{equation}\label{cotaabsoluta}
\left| \frac{Q_{n,r}^{(k)}(0)}{(L_n^{\alpha})^{(k)}(0)} \right|
\le 2(r+1), \quad  \,  \forall n \ge n_0  \quad  0\le  k \le n\,.
 \end{equation}
\end{lemma}
\textbf{Proof.} We prove the Lemma using mathematical  induction
on $k$, i.e., on the order of the derivative.

Keeping in mind Lemma \ref{c} (a) for $k=r+1$ and Lemma
\ref{lema2},  there exists a positive integer number $n_0$,
independent of $k$,  such that for all $n \ge n_0$, the following
formulas hold,

\begin{equation}\label{desig1}
 \left| \frac{Q_{n,r}^{(r+1)}(0)}{(L_n^{\alpha})^{(r+1)}(0)} \right|  \le 4 (r+1) (r+1)!  \frac{\Gamma(\alpha+r+2)}{\Gamma(\alpha+2r+3)} \,,
\end{equation}

\begin{equation}\label{desig2}
 \frac{\Gamma(\alpha+r+2)}{\Gamma(\alpha+i+1)} \left| \frac{A_{n,r}^{i}}{(L_n^{\alpha})^{(i)}(0)} \right|  \le1 \,, \quad i=1,\dots,r \,,
\end{equation}
and

\begin{equation}\label{desig3}
 \left| \frac{A_{n,r}^{r+1}}{(L_n^{\alpha})^{(r+1)}(0)} \right|  \le 2(r+1)!  \,.
\end{equation}

Notice that formula (\ref{desig1}) is  the required bound   for
$k=r+1$. Now, we assume that the result holds for a fixed $k$,
with $k \ge r+1 \,,$ and then  we will deduce that it holds for
$k+1$.

Taking absolute values in (\ref{ratioderiv})  and using induction hypothesis, (\ref{desig2}) and (\ref{desig3}), we get for $n \ge n_0$ and $k+1\le n$

\begin{align*}
& \left|\frac{Q_{n,r}^{(k+1)}(0)}{(L_n^{\alpha})^{(k+1)}(0)} \right| \le  \frac{\alpha+k+1}{\alpha+r+k+2}  \left| \frac{Q_{n,r}^{(k)}(0)}{(L_n^{\alpha})^{(k)}(0)} \right|    \\
& +    \frac{\Gamma(\alpha+k+2)}{\Gamma(\alpha+r+k+3)} \left[ \sum_{i=1}^r  \left( \begin{matrix} k \\ i-1  \end{matrix} \right) +  \left( \begin{matrix} k \\ r  \end{matrix} \right)  2(r+1)! \right] \\
&\le \left[ 2 (r+1) \frac{k!\,(k-(r-1))}{(k-r)!}  + 2 (r+1) \frac{k!}{(k-r)!} + \sum_{i=1}^r  \left( \begin{matrix} k \\ i-1  \end{matrix} \right)    \right] \frac{\Gamma(\alpha+k+2)}{\Gamma(\alpha+r+k+3)} \\
& \le \left[ 2 (r+1)\, \frac{k!}{(k-r)!}  \,(k+1-(r-1)) +r \frac{k!}{(k+1-r)!}  \right]  \frac{\Gamma(\alpha+k+2)}{\Gamma(\alpha+r+k+3)}  \\
& \le \left[ 2 (r+1)\, \frac{(k+1)!}{(k-r)!}   +\frac{(k+1)!}{(k+1-r)!}  \right]  \frac{\Gamma(\alpha+k+2)}{\Gamma(\alpha+r+k+3)}  \\
& \le 2 (r+1)\, \frac{(k+1)!}{(k+1-r)!}  \,(k+1-(r-1)) \frac{\Gamma(\alpha+k+2)}{\Gamma(\alpha+r+k+3)}  \,.
\end{align*}

So, the first part of Lemma is proved.  For the second part, using
(\ref{derivadaLn}), (\ref{Qn0(x)}) and  Lemma \ref{nucleos} (b) we
deduce for every $n\ge k$ the explicit expression

\begin{equation*}
\frac{Q_{n,0}^{(k)}(0)}{(L_n^{\alpha})^{(k)}(0)} = 1 - \frac{M_0
K_{n-1}(0,0)}{1+M_0 K_{n-1}(0,0)}
\frac{\alpha+1}{\alpha+k+1}\frac{n-k}{n} \,.
\end{equation*}
Then,
\begin{equation*}
0 < \frac{Q_{n,0}^{(k)}(0)}{(L_n^{\alpha})^{(k)}(0)} < 1 \,,
\end{equation*}
holds for all $k$ with $0\le k \le n$. According to this fact and
the first part of Lemma \ref{lema3}, we have (\ref{cotaabsoluta}).
 $\quad \Box $

 \bigskip

In the next theorem we show how the presence of the masses
  in  the inner product changes the asymptotic behaviour around  the
  origin. This result proves the conjecture posed in \cite{alv-mb}.

\begin{theorem}\label{teorMH}
 Let $\{ Q_{n,r}\}_{n \ge 0}$ be the SMOP with respect to the inner product (\ref{pr-ls1}). Then,
\begin{equation*}
\lim_n \frac{(-1)^n}{n!\,n^{\alpha}} Q_{n,r} \left(\frac{x}{n+j}\right) = (-1)^{r+1} x^{- \alpha/2} J_{\alpha+2r+2}(2 \sqrt{x}),
\end{equation*}
uniformly on compact  subsets of  $\mathbb{C}$ and uniformly on $j\in \mathbb{N}\cup \{0\}$.
\end{theorem}
\textbf{Proof.} Using the Taylor expansion of the polynomial $ Q_{n,r}$ we can write
\begin{equation*}
 \frac{(-1)^n}{n!\,n^{\alpha}} Q_{n,r}\left(\frac{x}{n+j}\right) = \sum_{k=0}^n \frac{(-1)^n}{n!\,n^{\alpha}}  \frac{Q_{n,r}^{(k)}(0)}{k!} \frac{x^k}{(n+j)^k} \,.
\end{equation*}
 To obtain the asymptotic behaviour of the above expression when $n\to \infty$,  we are going to use use the Lebesgue's dominated convergence theorem.
 Indeed, given a compact set $K \subset \mathbb{C}$, from
 (\ref{derivadaLn}), (\ref{equivgamma}), and   (\ref{cotaabsoluta}) in Lemma \ref{lema3}
  there exists a positive integer number $n_0$ such that for all $n \ge n_0$, for all $j \ge 0$ and for all $x\in K$,
\begin{align*}
& \frac{1}{n!\,n^{\alpha}} \left| \frac{Q_{n,r}^{(k)}(0)}{k!} \frac{x^k}{(n+j)^k} \right| \le  \frac{\Gamma(\alpha+1)|L_n^{\alpha}(0)|}{n!\,n^{\alpha}}  \frac{n!}{(n-k)!\,(n+j)^k} \frac{2(r+1)}{\Gamma(\alpha+k+1)} \frac{|x|^k}{k!}  \\
&\le \frac{4 (r+1)}{\Gamma(\alpha+k+1)} \frac{M^k}{k!}\,,
\end{align*}
\noindent for each $k\ge0$, where $M$ is a positive constant
depending on $K$. As $\displaystyle \sum_{k=0}^\infty
\frac{1}{\Gamma(\alpha+k+1)} \frac{M^k}{k!}$
 converges, the assumptions of the Lebesgue's dominated convergence theorem are satisfied.
 Then, using Lemma \ref{c} (a), (\ref {derivadaLn}), and (\ref{equivgamma}), we have
\begin{align*}
& \lim_n  \sum_{k=0}^n   \frac{(-1)^n}{n!\,n^{\alpha}} \frac{Q_{n,r}^{(k)}(0)}{k!} \frac{x^k}{(n+j)^k} \\ \nonumber
 & = \sum_{k=r+1}^{\infty} \frac{(-1)^k}{(k-(r+1))!} \frac{1}{\Gamma(\alpha+k+r+2)}  \, x^k
 = (-1)^{r+1} x^{- \alpha/2} J_{\alpha+2r+2}(2 \sqrt{x}) \,,
\end{align*}
uniformly  on  compact  subsets of $\mathbb{C}$ and uniformly on
$j\in \mathbb{N}\cup \{0\}.$ Thus,  the result follows. $ \quad
\Box $

\bigskip

In the next corollary we will show a remarkable difference between
the zeros of the orthogonal polynomials $L_n^{\alpha}$ and the
ones of  $Q_{n,r}$ concerning the convergence acceleration to $0$.

 Before analyzing this, we recall (see \cite{sz}) that the zeros of the Laguerre polynomials are real, simple
and they are located in $(0,\infty)$. We denote by
$(x_{n,k})_{k=1}^n$ the zeros of $L_n^{\alpha}$ in an increasing
order. It is worth pointing out that they satisfy the interlacing
property $0 < x_{n+1,1} < x_{n,1} < x_{n+1,2} < \dots $, and that
$x_{n,k} \underset{n} \to 0$ for each fixed $k$.

 Let $(j_{\alpha,k})_{k\ge 1}$ be the positive zeros of the Bessel function $J_{\alpha}$ in an increasing order.
Then, formula (\ref{MH Ln})  and Hurwitz's
theorem lead us to
 $$n x_{n,k}\underset{n} \to j_{\alpha,k}, \quad k\ge 1,$$
  and therefore
$$x_{n,k} \sim \frac{C_k}{n}, \quad k\ge 1,$$
where $C_k$ is a positive constant depending on $k$.

Concerning the zeros of $Q_{n,r}$,  standard arguments (see for
instance \cite{ch}) allow us to establish that $Q_{n,r}$ has at
least $n-(r+1)$ zeros with odd multiplicity  in $(0, +\infty)$, or
equivalently $n-(r+1)$ changes of sign. Moreover, since $M_0>0$
and the mass point in the discrete part of the inner product
belongs to the boundary of $(0, +\infty)$ then the number of zeros
with odd multiplicity is at least $n-r$ (see \cite{alr}).

From Theorem \ref{teorMH} and Hurwitz's theorem and taking into
account the multiplicity of 0 as a zero of the limit function in
Theorem \ref{teorMH} we achieve

\begin{corollary}\label{acelceros} Let $(\xi_{n,k}^{r})_{k=1}^n$ be the zeros of  $Q_{n,r}$. Then

 $$n\,\xi_{n,k}^{r} \underset{n} \to 0, \quad 1 \le k \le r+1,$$

 $$n\,\xi_{n,k}^{r} \underset{n} \to j_{\alpha+2r+2,k-r-1},  \quad k \ge  r+2.$$

\end{corollary}

\noindent\textbf{Remark 1.} The presence of the positive masses
$M_i, \, i=0, \ldots,  r,$ in the inner product produces a
convergence acceleration to $0$ of $r+1$ zeros of the polynomials
$Q_{n,r}$.

\subsection{Generalized Hermite--Sobolev polynomials}

As a consequence of the previous results, we are going to establish asymptotic  properties for the orthogonal  polynomials associated with  the following
inner product
\begin{equation}\label{pr-h}
(p,q) =\int_{\mathbb{R}}
p(x) q(x)\vert x \vert^{2\mu}\,e^{-x^2}dx+\sum_{i=0}^{2r+1} M_i\,p^{(i)}(0)\,q^{(i)}(0) ,
\end{equation}with $\mu>-1/2$ and $M_i > 0, \quad i=0, \dots, 2r+1$. We denote by $S_{n,r}^{\mu}$ their monic orthogonal polynomials.

The polynomials  $H_{n}^{\mu}$ orthogonal with respect to the
weight $\vert x \vert^{2\mu}\,e^{-x^2}$ ($\mu>-1/2$) are called
 \textit{generalized Hermite polynomials}.

Notice that in this case the polynomials $S_{n,r}^{\mu}$ are symmetric,
that is, $S_{n,r}^{\mu}(-x)=(-1)^n\,S_{n,r}^{\mu}(x),$ and because of this
symmetry, we can transform this inner product (\ref{pr-h}) into an
inner product like (\ref{pr-ls1}) and so we can establish a simple
relation between the polynomials $S_{n,r}^{\mu}$ and the polynomials
$Q_{n,r}$  considered before. This technique is known as a
symmetrization process. In fact, in \cite{ch} this process is
considered for standard inner products associated with positive
measures. The simplest case of this situation is the relation
between  monic Laguerre polynomials and Hermite polynomials, that
is (see \cite{ch} or \cite{sz}),

$$
H_{2n}(x)=L_n^{-1/2}(x^2), \quad H_{2n+1}(x)=xL_n^{1/2}(x^2),
\quad n\ge 0.
$$
 Later in \cite{ammr} the authors generalize the symmetrization process in the framework of Sobolev type
 orthogonal polynomials, see Theorem 2 in \cite{ammr}.

 As a consequence we have
 $$S_{2n,r}^{\mu}(x)=Q_{n,r}^{\mu-1/2}(x^2), \quad S_{2n+1,r}^{\mu}(x)=xQ_{n,r}^{\mu+1/2}(x^2)$$
 \noindent where $\{Q_{n,r}^{\mu-1/2}\}_{n\ge0}$ (respectively, $\{Q_{n,r}^{\mu+1/2}\}_{n\ge0}$) is the
 SMOP with respect to an inner product like (\ref{pr-ls1}) with $\alpha=\mu-1/2$ (respectively, $\alpha=\mu+1/2$).

 Thus, applying the above Theorem \ref{relative-asymptotic} and Theorem \ref{teorMH} in a straightforward way
 we obtain

\begin{proposition}\label{hermite}
 Let $\{ S_{n,r}^{\mu}\}_{n \ge 0}$ be the SMOP with respect to the inner product (\ref{pr-h}).
Then,
\begin{itemize}
\item[(a)]
\begin{equation*}
\lim_n \frac{S_{n,r}^{\mu}(x)}{H_n^{\mu}(x)} = 1,
\end{equation*}
uniformly  on  compact  subsets of $\mathbb{C}\setminus
\mathbb{R}$.

\item[(b)]
\begin{align*}
&\lim_n \frac{(-1)^n \sqrt{n}}{n!\,n^{\mu}} S_{2n,r}^{\mu} \left(\frac{x}{2 \sqrt{n+j}}\right) = (-1)^{r+1} \left( \frac{x}{2} \right)^{-\mu+1/2} J_{\mu+2r+3/2}(x), \\
&\lim_n \frac{(-1)^n }{n!\,n^{\mu}} S_{2n+1,r}^{\mu} \left(\frac{x}{2 \sqrt{n+j}}\right) = (-1)^{r+1} \left( \frac{x}{2} \right)^{-\mu+1/2} J_{\mu+2r+5/2}(x),
\end{align*}
\end{itemize}
uniformly  on  compact  subsets of $\mathbb{C}$ and uniformly on $j\in \mathbb{N}\cup \{0\}$.

\end{proposition}

\noindent\textbf{Remark 2.} These results generalize some of the
results in \cite{ampr} and solve the conjecture posed there.

\section{Inner products with holes}

In this section, we are concerned with inner products such that in their discrete
part  at least one of the masses vanishes:
\begin{equation}\label{pr-ls2}
(p,q)_{r,s} =(p,q)_r + M_s p^{(s)}(0)q^{(s)}(0), \quad s \ge r+2,
\end{equation}
where $M_s > 0$, and in $(. , . )_r $ all the masses are positive.
That is, roughly speaking, there is  a \lq\lq hole" in the
discrete part of  the inner product $(. , . )_{r,s}$. We denote by
$\{ T_{n,r,s} \}_{n\ge0}$  the sequence of monic polynomials
orthogonal with respect to the inner product $(. , . )_{r,s}$.

The Fourier expansion of the polynomial $T_{n,r,s}$ in the orthogonal basis $\{ Q_{n,r} \}_{n\ge0}$ gives
\begin{equation}\label{expanTns}
T_{n,r,s}(x) =Q_{n,r}(x) -  \frac{M_s Q_{n,r}^{(s)}(0)}{1+M_s K_{n-1,r}^{(s,s)}(0,0)}K_{n-1,r}^{(0,s)}(x,0) \,.
\end{equation}

\medskip
Using similar arguments as in Lemma \ref{c} it can be proved the following

\begin{lemma}\label{lema1aguj}
Let $\{ T_{n,r,s} \}_{n\ge0}$ be the SMOP with respect to the inner product  (\ref{pr-ls2}). Then
the following statements hold:

\begin{itemize}
\item[(a)]
For either $0 \le k \le r$ or $k=s$,
$$ \frac{T_{n,r,s}^{(k)}(0)}{(L_n^{\alpha})^{(k)}(0)} \sim  \frac{C_{r,s,k}}{n^{\alpha+2k+1}}\, ,$$
where $C_{r,s,k}$ is a nonzero real number independent of $n$.

 For $k \ge r+1$ and $k \not =s$
$$\lim_n \frac{T_{n,r,s}^{(k)}(0)}{(L_n^{\alpha})^{(k)}(0)}=
\frac{k!}{(k-(r+1))!} \frac{k-s}{\alpha +s+k+1}\frac{\Gamma(\alpha+k+1)}{\Gamma(\alpha+r+k+2)}  \, .$$

\item[(b)]
$$\lim_n \frac{(T_{n,r,s} T_{n,r,s})_{r,s}}{\Vert L_n^{\alpha} \Vert ^2} = 1 \,.$$
\end{itemize}
\end{lemma}

The above lemma also allows us to deduce the relative asymptotics
for these orthogonal polynomials.

\begin{theorem}

Let $\{ T_{n,r,s} \}_{n\ge0}$ be the SMOP with respect to the
inner product defined by (\ref{pr-ls2}). Then
\begin{equation*}
\lim_n \frac{T_{n,r,s}(x)}{L_n^{\alpha}(x)} = 1,
\end{equation*}
 uniformly on compact subsets of $\mathbb{C} \setminus [0,\infty)$.
\end{theorem}
\textbf{Proof.} Proceeding  as in Theorem 2 we have
\begin{equation*}
\frac{T_{n,r,s}(x)}{L_n^{\alpha}(x)} = 1 - \sum_{i=0}^r M_i
T_{n,r,s}^{(i)}(0) \frac{K_{n-1}^{(0,i)}(x,0)}{L_n^{\alpha}(x)} -
M_s T_{n,r,s}^{(s)}(0)
\frac{K_{n-1}^{(0,s)}(x,0)}{L_n^{\alpha}(x)}\, ,
\end{equation*}
and in this way, we can prove that each one of the terms in the
above sum converges to $0$ uniformly  on compact subsets of
$\mathbb{C} \setminus [0,\infty)$ and the result follows. $\Box$

The Mehler--Heine type formula for the polynomials $\{
T_{n,r,s}\}_{n \ge 0}$ is  also obtained now.

\begin{theorem}\label{teorMHaguj}
 Let $\{ T_{n,r,s}\}_{n \ge 0}$ be the SMOP with respect to the inner product (\ref{pr-ls2}). Then,
\begin{align}\label{limTn}
&\lim_n \frac{(-1)^n}{n!\,n^{\alpha}} T_{n,r,s}
\left(\frac{x}{n+j}\right)= (-1)^{r+1} x^{- \alpha/2}  \\
\nonumber & \times
\left[\frac{-(s-(r+1))}{\alpha+r+s+2}J_{\alpha+2r+2}(2 \sqrt{x}) +
\sum_{i=2}^{s-r+1} \lambda_i J_{\alpha+2r+2i}(2 \sqrt{x})  \right]
\,,
\end{align}
where $\lambda_i$ are nonzero real
numbers. The limit holds uniformly  on  compact  subsets of
$\mathbb{C}$ and uniformly on $j\in \mathbb{N}\cup \{0\}$.
\end{theorem}
\textbf{Proof.} From (\ref{expanTns}),

\begin{align}\label{Tn}
& \frac{(-1)^n}{n! \,n^{\alpha}} T_{n,r,s}\left(\frac{x}{n+j}\right) \\ \nonumber
& = \frac{(-1)^n}{n! \,n^{\alpha}} Q_{n,r}\left(\frac{x}{n+j}\right)
-  \frac{M_s Q_{n,r}^{(s)}(0)}{1+M_s K_{n-1,r}^{(s,s)}(0,0)} \frac{(-1)^n}{n! \,n^{\alpha}} \sum_{k=0}^{n-1} \frac{K_{n-1,r}^{(k,s)}(0,0)}{k!} \left(\frac{x}{n+j}\right)^k.
\end{align}

To estimate the kernels $K_{n-1,r}^{(k,s)}(0,0)$ we apply Stolz
criterion, Lemma \ref{c}, (\ref{derivadaLn}) and (\ref{estimate}),
obtaining

\begin{align*}
& \lim_n  \frac{K_{n-1,r}^{(k,s)}(0,0)}{n^{\alpha+k+s+1}} \\
&=\begin{cases}  0 & \text {if $0 \le k \le r$} \, ,\\
\frac{k!}{(k-(r+1))!}\frac{s!}{(s-(r+1))!}\frac{(-1)^{k+s}\Gamma(\alpha+1)}{(\alpha+k+s+1)\Gamma(\alpha+k+r+2)\Gamma(\alpha+s+r+2)} & \text {if $k \ge r+1$} \, .\end{cases}
\end{align*}

Moreover, it is not difficult to check that

\begin{equation*}
\lim_n  \frac{(-1)^{n} \, n^{s+1}}{n!}  \frac{Q_{n,r}^{(s)}(0)}{K_{n-1,r}^{(s,s)}(0,0)} = (-1)^s\,\frac{(s-(r+1))!}{s!} \frac{(\alpha+2s+1) \Gamma(\alpha+s+r+2)}{\Gamma(\alpha+1)}.
\end{equation*}

According to  the two above results, we get the asymptotic
behaviour of the coefficients in the sum appearing in (\ref{Tn}),
\begin{align*}
&\lim_n \frac{(-1)^n}{n! \,n^{\alpha+k}} \frac{M_s Q_{n,r}^{(s)}(0) K_{n-1,r}^{(k,s)}(0,0)}{1+M_s K_{n-1,r}^{(s,s)}(0,0)} \\
&=\begin{cases}  0 & \text {if $0 \le k \le r$} \, ,\\
\frac{(-1)^{k}k!}{(k-(r+1))!} \frac{(\alpha+2s+1)}{(\alpha+k+s+1) \Gamma(\alpha+k+r+2)} & \text {if $k \ge r+1$} \, .
\end{cases}
\end{align*}

On the other hand, from (\ref{cotaabsoluta}) there exists $n_0 \in
\mathbb{N}$ such that  for all $n \ge n_0$, and for every $ k,s \ge0 $ we
have

\begin{equation*}
\left|K_{n-1,r}^{(k,s)}(0,0) \right| \le 4 (r+1)^2 \left|K_{n-1}^{(k,s)}(0,0) \right| \,.
\end{equation*}

Now, to obtain a bound for the kernels $K_{n-1}^{(k,s)}(0,0)$, we
consider  the expression which appears  in Lemma \ref{nucleos}
(b). Then, when $k \ge s-1$, it is easy to check that $j! \,
(k+s+1-j)! \ge s! \, (k+1)!$, and $\Gamma(\alpha+j+1) \,
\Gamma(\alpha+k+s+2-j) \ge \Gamma(\alpha+s+1) \,
\Gamma(\alpha+k+2)$ for $0 \le j \le s$. Therefore,
\begin{equation*}
\left| K_{n-1}^{(k,s)}(0,0) \right| \le C_s\, \frac{n^{\alpha+k+s+1}}{\Gamma(\alpha+k+2)}  \,.
\end{equation*}

Indeed, given a compact set $K \subset \mathbb{C}$, from
 (\ref{derivadaLn}),  (\ref{equivgamma}), and
 (\ref{cotaabsoluta}),
there exists a positive integer number $n_0$ such that for all $n
\ge n_0$, for all $k \ge s-1$, for all $j\ge0$ and for all $x\in
K$,
\begin{align*}
& \frac{1}{n^{\alpha +s+1}} \left| \frac{K_{n-1,r}^{(k,s)}(0,0)}{k!} \frac{x^k}{(n+j)^k} \right| \le  C_{s}\frac{4 (r+1)^2}{\Gamma(\alpha+k+2)} \frac{M^k}{k!},
\end{align*}
\noindent where $M$ is a positive constant depending on $K$. As $
\sum_{k=0}^\infty  \frac{1}{\Gamma(\alpha+k+2)} \frac{M^k}{k!}$
converges,  we can apply of the Lebesgue's dominated convergence
in the last term of (\ref{Tn}). Then, using  Theorem \ref{teorMH},
we obtain
\begin{align}\label{clBessel}
& \lim_n \frac{(-1)^n}{n! \,n^{\alpha}} T_{n,r,s}\left(\frac{x}{n+j}\right)  = (-1)^{r+1}x^{-\alpha/2} J_{\alpha+2r+2}(2\sqrt{x}) \\ \nonumber
& - \sum_{k=r+1}^{\infty} \frac{(-1)^{k}  (\alpha+2s+1)}{(\alpha+k+s+1) \Gamma(\alpha+k+r+2)} \frac{x^{k}}{(k-(r+1))!} .
\end{align}
If we write $s=r+1+h$ with $h\ge1$, then the above series is read
as

\begin{equation*}
(-1)^{r+1}x^{r+1}(\alpha+2r+2h+3)\sum_{k=0}^{\infty}
\frac{\Gamma(\alpha+k+2r+h+3)}{\Gamma(\alpha+k+2r+3)\Gamma(\alpha+k+2r+h+4)}
\frac{ (-1)^k x^k}{k!}.
\end{equation*}
 Observe that $
\Gamma(\alpha+k+2r+h+3)/\Gamma(\alpha+k+2r+3)$ is a polynomial in
$k$ of degree $h$ (the number of holes) and so we can write
\begin{equation*}
\frac{\Gamma(\alpha+k+2r+h+3)}{\Gamma(\alpha+k+2r+3)}=\frac{\Gamma(\alpha+2r+h+3)}{\Gamma(\alpha+2r+3)}
+\sum_{l=1}^h \beta_{l}k^{l}
\end{equation*}
where $\beta_{l}, \, l=1,\dots,h$ are positive coefficients.
Thus, the above series can be expressed as
\begin{align*}
&\frac{\Gamma(\alpha+2r+h+3)}{\Gamma(\alpha+2r+3)}\sum_{k=0}^{\infty}
\frac{1}{\Gamma(\alpha+k+2r+h+4)} \frac{ (-1)^k x^k}{k!} \\&+
\sum_{l=1}^{h}\beta_{l}\sum_{k=0}^{\infty}\frac{k^{l}}{\Gamma(\alpha+k+2r+h+4)}\frac{
(-1)^k  x^k}{k!}.
\end{align*}
For the first one, using  the recurrence relation repeatedly (see
\cite{sz}),
\begin{equation*}\label{Bessel}
J_{\alpha -1}(2 \sqrt{x})+J_{\alpha +1}(2
\sqrt{x})=\alpha\,x^{-\frac{1}{2}}\,J_{\alpha}(2 \sqrt{x}),
\end{equation*}
we get
\begin{align*}
&\frac{\Gamma(\alpha+2r+h+3)}{\Gamma(\alpha+2r+3)} \sum_{k=0}^{\infty}
\frac{1}{\Gamma(\alpha+k+2r+h+4)} \frac{ (-1)^k x^k}{k!}  \\ \nonumber
&=\frac{\Gamma(\alpha+2r+h+3)}{\Gamma(\alpha+2r+3)} \,
x^{-\frac{\alpha+2r+h+3}{2}} \, J_{\alpha+2r+h+3}(2\sqrt{x}) \\
\nonumber &= x^{-\frac{\alpha+2r+2}{2}} \left[
\frac{1}{\alpha+2r+h+3} J_{\alpha+2r+2}(2\sqrt{x}) +
\sum_{i=2}^{h+2} \mu_{i} J_{\alpha+2r+2i}(2\sqrt{x}) \right]\,,
\end{align*}
where $\mu_{i}$ are real numbers which can be computed explicitly.

Moreover, for the remaining series, using the same arguments it
can be seen that each one of the terms can be written as a
combination of Bessel functions
 of order bigger than $\alpha+2r+2$. More precisely, for $l=1, \dots, h$
\begin{align*}
& \sum_{k=0}^{\infty} \frac{k^l}{\Gamma(\alpha+k+2r+h+4)} \frac{
(-1)^k  x^k}{k!} =-x \sum_{k=0}^{\infty}
\frac{(k+1)^{l-1}}{\Gamma(\alpha+k+2r+h+5)} \frac{ (-1)^k x^k}{k!}
\\ \nonumber &= x^{-\frac{\alpha+2r+2}{2}}\left[
-\frac{\Gamma(\alpha+2r+5)}{\Gamma(\alpha+2r+h+5)}J_{\alpha+2r+4}(2\sqrt{x})
+ \sum_{i=3}^{h+2} \mu^{*}_{i} J_{\alpha+2r+2i}(2\sqrt{x})
\right],
\end{align*}
where  $\mu_{i}^*$ are again real numbers which can be computed
explicitly.

Finally, taking these results into account in (\ref{clBessel}) we
achieve

\begin{align*}
 \lim_n \frac{(-1)^n}{n! \,n^{\alpha}} T_{n,r,s}\left(\frac{x}{n+j}\right) & = (-1)^{r+1}x^{-\alpha/2}\left[1-\frac{\alpha+2r+2h+3}{\alpha+2r+h+3}\right] J_{\alpha+2r+2}(2\sqrt{x}) \\
&+(-1)^{r+1}x^{-\alpha/2}\sum_{i=2}^{h+2} \lambda_i J_{\alpha+2r+2i}(2\sqrt{x}),\\
\end{align*}
and the proof is concluded. $\quad \Box$

\bigskip

We can also use the techniques developed in this section  to
obtain asymptotics results in other similar frameworks. For
example, very recently in \cite{dm1} the authors consider the
inner product
\begin{align} \label{herpac}
(p,q)_* &=\frac{1}{\Gamma(\alpha +1)}\int_0^{\infty } p(x) q(x)\,
x^{\alpha} e^{-x} \,dx + M_s p^{(s)}(0)q^{(s)}(0)\nonumber \\
&=(p,q)+M_s p^{(s)}(0)q^{(s)}(0),
\end{align}
where $s\ge 1.$ Notice that the role played by $(p,q)_r$ in the
inner product (\ref{pr-ls2}) is now played by $(p,q)$ in
(\ref{herpac}). Therefore, if we  proceed as in this section, we
improve  the asymptotic results appearing in \cite{dm1}.

\medskip

Moreover, for the particular case $s=r+2\, $ in the inner
product (\ref{pr-ls2}), i.e., when there is a hole of ``length
one'', the result established in the above theorem generalize the
one obtained in \cite{alv-mb}. In fact, handling the right--hand
side of the expression (\ref{clBessel}) and using the recurrence
relation of the Bessel functions we obtain:

\begin{align*}
& \lim_n \frac{(-1)^n}{n! \,n^{\alpha}} T_{n,r,r+2}\left(\frac{x}{n+j}\right)= (-1)^{r+1} x^{- \alpha/2} \\
&\times
\left[\frac{-1}{\alpha+2r+4}J_{\alpha+2r+2}(2 \sqrt{x})
-J_{\alpha+2r+4}(2 \sqrt{x})+ \frac{-1}{\alpha+2r+4}J_{\alpha+2r+6}(2 \sqrt{x}) \right].
\end{align*}

\medskip
On the other hand, as a consequence of the above theorem we
present the situation about the acceleration of the convergence
towards the origin of the zeros of the polynomials $T_{n,r,s}$.
The quasi--orthogonality of order $s+1$ of the sequence
$\{T_{n,r,s}\}_{n\ge0}$ with respect to the positive measure
$x^{\alpha+s+1}e^{-x}$ assures that $T_{n,r,s}$ has at least
$n-(s+1)$ changes of sign in $(0,+\infty)$. However, in \cite{alr}
the authors proved that the number of zeros in $(0, +\infty)$ does
not depend on the order of the derivatives but on the number of
terms in the discrete part of the inner product. So, $T_{n,r,s}$
has at least $n-(r+1)$ zeros with odd multiplicity in
$(0,+\infty)$.

From Theorem \ref{teorMHaguj} and Hurwitz's theorem and taking
into account that $x=0$ is a zero of multiplicity $r+1$ of the
limit function in (\ref{limTn}) we achieve the following result:

\begin{corollary} Let $(\zeta_{n,k}^{r,s})_{k=1}^n$ be the zeros of  $T_{n,r,s}$. Then

$$ n\,\zeta_{n,k}^{r,s} \underset{n} \to 0, \quad 1 \le k \le r+1,$$

$$ n\,\zeta_{n,k}^{r,s} \underset{n} \to j_{\alpha+2r+2,k-r-1}, \quad k \ge r+2.$$

\end{corollary}

\medskip

\noindent  In the next remark, we compare  the previous result
with the corresponding one of Corollary \ref{acelceros}.

\noindent\textbf{Remark 3.}  We want to highlight that the
convergence acceleration to $0$ of the zeros of the polynomials
$Q_{n,r}$  and  $T_{n,r,s}$ is the same. That is, the addition of
a mass $M_s$ \textit{after a hole} in the inner product does not
affect  the convergence acceleration to $0$.

\bigskip

As we have explained  in the previous section, using a
symmetrization process, we can obtain the relative asymptotics and
the Mehler--Heine type formulas for generalized Hermite--Sobolev
polynomials with holes in the discrete part of the  inner product.

\end{document}